\documentclass[a4paper,12pt]{article}
\usepackage{amsmath,amsfonts}
\date{16 june 2002}
\title{Embeddings of ultradistributions and periodic hyperfunctions
in Colombeau type algebras through sequence spaces}

\author{Antoine Delcroix\\\small\tt Antoine.Delcroix@univ-ag.fr
\and Maximilian F. Hasler\thanks{
corresponding author:\newline Univ. Antilles--Guyane, D.S.I.,
B.P. 7209, 97275 Schoelcher cedex (Martinique, F.W.I.)
}\\\small\tt MHasler@martinique.univ-ag.fr
\and Stevan Pilipovi\'c\\\small\tt pilipovic@im.ns.ac.yu
\and Vincent Valmorin\\\small\tt Vincent.Valmorin@univ-ag.fr}

\makeatletter
\advance\textheight\topmargin
\newtheorem{theorem}{Theorem}[section]

\newtheorem{definition}[theorem]{Definition}
\newtheorem{proposition}[theorem]{Proposition}
\newtheorem{lemma}[theorem]{Lemma}
\newtheorem{example}{Example}
\newtheorem{remark}{Remark}
\newenvironment{proof}[1][Proof]{\par\nb{#1.} }{\hfill$\Box$\par\vskip1ex}
\newcommand\PROOF[1]{\begin{proof}#1\end{proof}}

\newcommand\A{{\cal A}}
\newcommand\B{{\cal B}}
\newcommand\CC{{\cal C}}
\newcommand\C{{\mathbb C}}
\newcommand\comment[1]{}
\newcommand\CEP{\ensuremath{(\mathcal C,\mathcal E,\mathcal P)}}
\newcommand\D{{\cal D}}
\newcommand\DD{{\textbf{\textsl D}}}
\newcommand\der{{\rm der}}
\newcommand\E{{\cal E}}
\newcommand\e{\varepsilon}
\newcommand\eN{{\cal N}}
\newcommand\F{{\cal F}}
\newcommand\FT[1]{\mathop{\mathcal{F\!T}}(#1)}
\newcommand\G{{\cal G}}

\def\INT#1d{\int#1\,\mathrm d}
\newcommand\J{{\cal J}}
\newcommand\ind{\mathop{\rm ind\,lim}}
\newcommand\K{{\cal K}}
\renewcommand\l{\lambda}
\newcommand\lr[3]{\left#1#3\right#2}
\newcommand\Mid{~\Big|~}
\newcommand\N{{\mathbb N}}
\newcommand\Ns{{\mathbb N^*}}
\newcommand\nb[1]{\noindent\textbf{#1}}
\renewcommand\O{{\cal O}}
\newcommand\ola[1]{\rlap{$\displaystyle\overleftarrow{\phantom{#1}}$}#1}
\newcommand\ora[1]{\rlap{$\displaystyle\overrightarrow{\phantom{#1}}$}#1}
\newcommand\olra{\overleftrightarrow}
\renewcommand\P{{\cal P}}
\newcommand\p[1]{\left(#1\right)}
\newcommand\pow{{\rm pow}}
\newcommand\R{{\mathbb R}}
\newcommand\real{\mathop{\Re\!e}}
\newcommand\supp{\operatorname{supp}}
\newcommand\set[1]{\ensuremath{\left\{\,#1\,\right\}}}
\newcommand\T{{\mathbb T}}
\newcommand\ultra[1]{|\!|\!|\,#1\,|\!|\!|}
\newcommand\X{{\cal X}}
\newcommand\Z{{\mathbb Z}}

\begin{document}
\maketitle
\begin{abstract}
In a recent paper
, we gave
a topological description of Colombeau type
algebras introducing algebras of sequences with exponential weights.
Embeddings of Schwartz' spaces into Colombeau algebra $\G$ are well
known, but for ultradistribution and periodic hyperfunction type
spaces we give new constructions. We show that the multiplication of
regular enough functions (smooth, ultradifferentiable or
quasianalytic), embedded into corresponding algebras, is the ordinary
multiplication.
%
\\[1ex]
\end{abstract}

\section{Introduction}

Differential algebras of generalized functions containing embedded 
distributions are a convenient framework for the analysis of problems 
with singular coefficients and/or singular data, especially for the
non linear problems since the multiplications and other non linear 
operations are in general not defined in classical generalized function
spaces. 
Nowdays, there is a considerable litterature concerning such algebras. 
(For example see~\cite{bia,co,co1,gros,nepisca,ro2} and the references 
therein.) 

We have proved in~\cite{dhpv1} that those algebras, refered here as Colombeau
type algebras, can be reconsidered as a class of sequence spaces algebras 
and we gave a purely topological description of them.

In analogy  to embeddings of Schwarz' distributions, we show in this paper 
that some classes of ultradistributions and periodic hyperfunctions 
can be embedded into well choosen sequence algebras. 
Morever, we show that the product of enough regular elements,
ultradifferentiable functions or quasianalytic periodic functions of
appropriate classes, is the ordinary multiplication.

The problem of embedding of classical spaces into
corresponding sequence spaces algebra is closey related to the choice of
sequences of mollifiers, sequences of appropriately smooth functions 
converging to the delta distribution. While such a problem is trivial 
for embedding of Schwartz distributions it is essential for 
ultradifferentiable functions and ultra
distributions, considered in section 3. 
The same holds for periodic quasianalytic functions and correspondaing
periodic hyperfunctions of section 4.


Colombeau ultradistributions corresponding to a general non-quasianalytic
sequence were introduced and analyzed in \cite{pisca}.
Although we consider here the Gevrey sequence ${(p!^m)}_p$, $m>1$,
we give sharper estimates and improve results of \cite{pisca}:
%
The construction of appropriate mollifiers enable us to
give more precise results concerning embeddings.\comment{improve...}
Colombeau periodic hyperfunctions introduced in this paper are more
closely related to the global theory of generalized functions than
those of~\cite{valm}.
In this sense, we improve results of~\cite{valm}.


The novelty of results related to both cited papers and the embedding of both
classes of algebras into corresponding sequence space algebras,
and so their topological description, are the main results of this paper.

\section{General construction \cite{dhpv1}} \label{sect-GC}
In the sequel, we will use the following notations: $\N = \set{0,1,...}$ and
$\Ns = \set{1,2,...}$.

We recall our construction from \cite{dhpv1} for the case of $E=\C$,
$r_n=1/n^m, n\in\Ns$, where $m > 0$ is fixed. 
Let
$$
        \E^{m}_{0} = \set{c=(c_n)_n \in \C^{\Ns}~,~~\ultra{c}_{|\cdot|,1/n^m} =
  \limsup_{n\to\infty} |c_n|^{1/n^m} < \infty~}~.
$$
$$
{\eN}^{m}_{0} = \set{c=(c_n)_n \in \C^{\Ns}~,~~\ultra{c}_{|\cdot|,1/n^m} = 0~}~.
$$
The factor algebra ${\overline{\C}}^{m}=\E^{m}_{0}/{\eN}^{m}_{0}$,
$m > 0$, is
called the ring of Colombeau ultracomplex numbers for $m >1$ and the ring of
Colombeau hypercomplex numbers for $m \leq 1$.

Now we come to the general construction. 
Let $\left( E_\nu^\mu,p_{\nu\,}^\mu\right)_{\mu,\nu\in\Ns}$ be
a family of semi-normed algebras over $\R$ or $\C$ such that
\begin{align*}
        \forall\mu,\nu\in\Ns :~ E_\nu^{\mu+1} \hookrightarrow E_\nu^\mu
\text{ ~ and ~ }
        E_{\nu+1}^\mu \hookrightarrow E_\nu^\mu
\text{~ (resp. }
        E_\nu^\mu \hookrightarrow E_{\nu+1}^\mu ~)~,
\end{align*}
where $\hookrightarrow$ means continuously embedded.
~(For the $\nu$ index we consider inclusions in the two directions.)~
%
%
Then let
$\displaystyle 
        \ola{E}
        =\projlim_{\mu\to\infty} \projlim_{\nu\to\infty} E_\nu^\mu
        =\projlim_{\nu\to\infty} E_\nu^\nu
$, (resp.
$\displaystyle
        \ora{E}
        =\projlim_{\mu\to\infty} \ind_{\nu\to\infty}E_\nu^\mu)
$.
Such projective and inductive limits are usually considered with norms
instead of seminorms, and with the additional assumption that in the
projective case sequences are reduced, while in the inductive case for
every $\mu\in\Ns$ the inductive limit is regular, i.e. a set
$A\subset{\ind\limits_{\nu\to\infty}}E_\nu^\mu$ is bounded iff
it is contained in some $E^\mu_\nu$ and bounded there.

Let  consider a positive sequence $ r={(r_n)}_n\in(\R_+)^\Ns$ ($\R_+= (0,+\infty)$) decreasing to
zero and define (with $p\equiv\p{p_\nu^\mu}_{\nu,\mu}$)
\begin{align*}
        \ola\F_{p,r} &= \set{ f\in\ola{E}^\Ns \Mid
        \forall\mu,\nu\in\Ns : \ultra f_{p_\nu^\mu,\,r} < \infty } ~,
\\
        \ola\K_{p,r} &= \set{ f\in\ola{E}^\Ns \Mid
        \forall\mu,\nu\in\Ns : \ultra f_{p_\nu^\mu,\,r} = 0 }
\\
\text{(resp.}~~
        \ora\F_{p,r} &= \bigcap_{\mu\in\Ns} \ora\F_{p,r}^\mu
~,~~    \ora\F_{p,r}^\mu = \bigcup_{\nu\in\Ns}
        \set{ f\in\p{ E_\nu^\mu }^\Ns \Mid \ultra f_{p_\nu^\mu,r} < \infty }
~,\\
        \ora\K_{p,r} &= \bigcap_{\mu\in\Ns} \ora\K_{p,r}^\mu
~,~~    \ora\K_{p,r}^\mu = \bigcup_{\nu\in\Ns}
        \set{ f\in\p{ E_\nu^\mu }^\Ns \Mid \ultra f_{p_\nu^\mu,r} = 0 }
~\text)~.
\end{align*}
Recall~\cite{dhpv1}:

{\def\item[#1]{}\parindent0pt
\item[(i)] Writing $\olra\cdot$ for both, $\ola{~\cdot}~$ or
$\ora{~\cdot}~$, we have that $\olra\F_{p,r}$ is an algebra and
$\olra\K_{p,r}$ is an ideal of $\olra\F_{p,r}$; thus, $\olra\G_{p,r}=
\olra\F_{p,r}/\olra\K_{p,r}$ is an algebra.

\item[(ii)] For every $\mu,\nu\in\Ns,~ d_{p_\nu^\mu}: {(E_\nu^\mu)}^\Ns
\times {(E_\nu^\mu)}^\Ns \to\overline\R_+$ defined by
$d_{p_\nu^\mu}(f,g)= \ultra{f-g}_{p_\nu^\mu,r}$ is an
ultrapseudometric on ${(E_\nu^\mu)}^\Ns$.  Moreover,
$(d_{p_\nu^\mu})_{\mu,\nu}$ induces a topological algebra\footnote
{\label{a} over $(\C^\Ns,\ultra\cdot_{|\cdot|})$, not over $\C$:
scalar multiplication is not continuous.}
structure on $\ola\F_{p,r}$
such that the intersection of
the neighborhoods of zero equals $\ola\K_{p,r}$.

\item[(iii)] From the properties above, the factor space $\ola\G_{p,r}=
\ola\F_{p,r}/\ola\K_{p,r}$ is a topological algebra over generalized
numbers $\overline{\C}_r=\G_{|\cdot|,r}$ (constructed with the sequence $r=r_n$
as above for the Colombeau ultracomplex numbers). The topology of
$\ola\G_{p,r}$ is defined by the family of
ultrametrics $(\tilde{d}_{p_\nu^\mu})_{\mu,\nu}$ where
$\tilde{d}_{p_\nu^\mu} ([f],[g])=d_{p_\nu^\mu}(f,g)$, $[f]$ standing
for the class of~$f$.

\item[(iv)] If $\tau_\mu$ denotes the inductive limit topology on
$\F_{p,r}^\mu=\bigcup_{\nu\in\Ns}((\tilde{E_\nu^\mu})^\Ns,d_{\mu,\nu})$,
$\mu\in\Ns$, then $\ora\F_{p,r}$ is a topological algebra
%
(over $(\C^\Ns,\ultra\cdot_{|\cdot|})$, as for $\ola\F_{p,r}$)
for the
projective limit topology of the family $(\F_{p,r}^\mu,\tau_\mu)_\mu$.

}

\begin{remark}
The two multiplicative sets $H=[0,1]$ and $I=[0,1)$ verify the relations
$H\cdot H = H,~ I\cdot H = I,~ I\cdot I = I$, just like the sets
$[0,\infty)$ and $\set0$.
Thus, similar constructions can also be made with $\ultra\cdot\le1$
and $\ultra\cdot<1$ instead of $\ultra\cdot<\infty$ and $\ultra\cdot=0$.
This is used in the setting of infra-exponential algebras and
also appears in the
context of periodic hyperfunctions. 
\end{remark}


\section{Colombeau ultradistributions of Gevrey class}
\label{kopi-ref}\label{Sect-UD}
\subsection{Ultradistributions of Gevrey class}

We refer to \cite{kom1} for definitions of the spaces
$\E^{(m)}$,$\D^{(m)}$, $\E^{\{m\}}$, $\D^{\{m\}}$ ($m>1$), and their
duals, Beurling and Roumieu type ultradistribution spaces.
Here we construct Colombeau ultradistribution algebras corresponding
to $M_p=p!^m$, $m>1$.  We apply the construction of section~\ref{sect-GC}.

For the function space $E=\CC^\infty(\R^s)$, we define
for all $\mu,\nu\in\R_+$ and $m>1$ the seminorms
$$
        p_\nu^{m,\mu}(f) = \sup_{ |x|\leq\mu , \alpha\in\N ^s}
        \dfrac{\nu^{|\alpha|}}{{\alpha!}^m} |f^{(\alpha)}(x)|
\text{~~ and ~~}
        q^{m,\mu}_{\nu}=p^{m,\mu}_{1/\nu} ~.
$$
Then let, for $\mu,\nu\in\Ns$, $E_\nu^\mu = E_{p^{m,\mu}_{\nu}}$ (resp.
$E_\nu^\mu = E_{q^{m,\mu}_{\nu}}$) be the subset of $ E $ on which
the given seminorm is finite.
For the first case, we clearly have $E_{\nu}^{\mu+1}\hookrightarrow
E_\nu^\mu$, $E_{\nu+1}^{\mu}\hookrightarrow E_\nu^\mu$
and for the second case, we have $E_{\nu}^{\mu+1}
\hookrightarrow E_{\nu}^{\mu}$, $E_{\nu}^{\mu}\hookrightarrow
E_{\nu+1}^{\mu}$ for any $\mu,\nu\in\Ns$.
Let $m>1$,  $m'>0$ and  $r_n=n^{-1/m'}$.\comment
{$m'-1\to m'$ but $m'>1$ unchanged}
\begin{definition}\label{D11}
The sets of exponentially growth order ultradistribution nets and null nets
of Beurling type are defined, respectively, by
$$
        \E^{(p!^m,p!^{m'})}_{exp} = \ola\F_{p^m,r}~,~~    
        \eN^{(p!^m,p!^{m'})} = \ola\K_{p^m,r} ~.
$$
The sets of exponentially growth order ultradistribution nets and null nets of
Roumieu type are defined, respectively, by
$$
        \E^{\{p!^m,p!^{m'}\}}_{exp} = \ora\F_{q^m,r}~,~~    
        \eN^{\{p!^m,p!^{m'}\}} = \ora\K_{q^m,r} ~.
$$
\end{definition}

\begin{proposition}
\begin{enumerate}
\item[(i)] $\E^{(p!^m,p!^{m'})}_{exp}$ (resp. $\E^{\{p!^m,p!^{m'}\}}_{exp}$)
are algebras under pointwise multiplication, and $\eN^{(p!^m,p!^{m'})}$
(resp.  $\eN^{\{p!^m,p!^{m'}\}}$) are ideals of them.
\item[(ii)] The pseudodistances induced by $\ultra\cdot_{p^{m,\mu}_\nu,m'}$
(resp.  $\ultra\cdot_{q^{m,\mu}_{\nu},m'}$) are ultrapseudometrics on
respective domains.
\end{enumerate}
\end{proposition}

\PROOF{With Definition~\ref{D11}, this is just a particular case of
of the general construction recalled in section~\ref{sect-GC}.}

The Colombeau ultradistribution algebra $\G^{(p!^m,p!^{m'})}$
(resp. $\G^{\{p!^m,p!^{m'}\}}$) is defined by
\begin{align*}
  \G^{(p!^m,p!^{m'})} = \ola\G_{p,r} 
  &= \E^{(p!^m,p!^{m'})}_{exp}/\eN^{(p!^m,p!^{m'})}
\\
\text{(resp.~ } \G^{\{p!^m,p!^{m'}\}} = \ora\G_{p,r}
  &=\E^{\{p!^m,p!^{m'}\}}_{exp}/\eN^{\{p!^m,p!^{m'}\}} ~\text)~.
\end{align*}
These topological algebras are invariant under the actions of
ultradifferential operators of respective classes $(m)$
and $\{m\}$, see e.g.~\cite{kom1}.

\subsection{Embeddings of ultradifferentiable functions and
  ultradistributions.}

%
In what follows, mollifiers will be constructed by elements of spaces
$\Sigma^\pow$ and $\Sigma_\der$, which consist of smooth functions
$\varphi$ on $\R$ with the property that for some $b>0$,
\[
  \sigma^b(\varphi) = \sup_{\beta\in\N, x\in\R} 
  \frac{ |x^\beta\,\varphi(x)| }{ b^\beta\,\beta! }  < \infty
\text{ ~ (resp. ~}
  \sigma_b(\varphi) = \sup_{\alpha\in\N, x\in\R}
  \frac{ |\varphi^{(\alpha)}(x)| }{ b^\alpha\,\alpha! } < \infty ~)~.
\]
Both spaces are endowed with the respective inductive topologies.
\begin{definition}
Let $\p{\phi^n}_{n\in\N^*}$ be a bounded net in $\Sigma^\pow$
(resp. $\Sigma_\der$) such that
$\forall n\in\N^*:\INT_\R t^j\phi^n(t)dt = \delta_{j,0}$ 
for $j\in\big\{ 0,1,2,\dots, [n^{1/m}]+1\big\},~ m>1$.
Then ${(\phi_n)}_{n\in\N^*}$ with $\phi_n = n\,\phi^n(n\,\cdot)$
is called a net of $\{m,\pow\}$ (resp. $\{m,\der\}$)--mollifiers.
\end{definition}
The following important lemma gives an explicit net of
$\{m,\pow\}$-- and $\{m,\der\}$--mollifiers:
\begin{lemma}
For all $n\in \N^*$ and $x\in \R$, let
$$
        h_n(x)=\exp\p{n^2-\sqrt[n]{n^{2n}+x^{2n}}} ~,~~
        k_n(x)=\exp\p{-x^{2n}} ~.
$$
Then, for all $n\in \N^*$, $h_n(0)=k_n(0)=1$ and 
$$
        \forall
        \alpha\in\set{1,...,2n-1}:~ 
        h_n^{(\alpha)}(0) = k_n^{(\alpha)}(0)= 0 ~,
$$
and there exist $r>0$ and $C>0$ such that
\begin{equation}
        \sup_{n\in\N^*} \sigma_r(h_n) < C ~,~~
        \sup_{n\in\N^*} \sigma^r(k_n) < C ~.\label{de,mu}
\end{equation}
Moreover, for given $m>1$, the nets\footnote
{we denote the Fourier transform by $\FT{\cdot}$ to avoid confusion with
spaces $\F_{p,r}$ etc.}
$$
        \phi^n = \frac1{2\,\pi}\FT{h_{g(n)}} 
\text{~ and ~}
        \phi^n = \frac1{2\,\pi}\FT{k_{g(n)}} ~,
$$
where $g(n)=\frac12[n^{1/(m-1)}]+1$ for $n\in\N^*$,
define a net of $\{m,\pow\}$--mollifiers and a net of
$\{m,\der\}$--mollifiers, respectively.
\end{lemma}

\begin{proof}
The first claims, $h_n^{(\alpha)}(0)=k_n^{(\alpha)}(0)=\delta_{\alpha,0}$ 
are easily verified, and imply obviously
$\int x^p \FT{h_n}=\int x^p \FT{k_n}=2\pi\,\delta_{p,0}~~ 
\forall p\in\set{0,...,2n-1}
$, which gives the second 
condition on \set{m,\der} resp.  \set{m,\pow}--mollifiers for $\phi^n$.

So let us show (\ref{de,mu}), i.e. $h_n\in\Sigma_\der$, $k_n\in\Sigma^\pow$
with constants independent of $n$. Consider first $h_n$.

The function $\C\ni z\mapsto\sqrt[n]{n^{2n}+z^{2n}}$ has 
singularities at $z=n\,e^{i\,\pi\,(2k+1)/(2n)}$. 
The nearest one to the real axis has the imaginary part
$n\sin\frac\pi{2n}\ge1~~\forall n\in\N^*$.
So for every $x\in\R$, the open disc $\set{|z-x|<1}$ lies in the domain
of analyticity of $h_n$.
Applying Cauchy's integral formula, we have
\begin{eqnarray*}
  \forall x\in\R, \forall n\in\N^* :~
  |h_n^{(\alpha)}(x)|  &=&  \left| \frac{\alpha!}{2\pi i}
    \int
_{|\zeta-x|=\frac12}
    \frac{ h_n(\zeta)\,\mathrm d\zeta }{ (\zeta-x)^{\alpha+1} } \right|
\\[1ex]
  &\le& 2^\alpha \, \alpha! \, \max_{\theta\in [0,2\pi]}
    \left| h_n( x+\tfrac12 e^{i \theta} ) \right|
~.
\end{eqnarray*}
Thus we have $\sigma_2(h_n)\le C$ and therefore (\ref{de,mu}), if
$\max|h_n(x+\frac12e^{i\theta})|<C$.
So let us show that there exists $C>0$ such that
\begin{equation}\label{br}
  \forall n\in\N^*, x\in\R:
     \real\p{n^2-\sqrt[n]{n^{2n}+(x+\tfrac12\,e^{i\theta})^{2n}}} < \ln C ~.
\end{equation}
Let $x+\frac12\,e^{i\,\theta}=\rho\,e^{i\,\phi}$
with $\rho\in\R,~ |\phi|<\frac\pi2$.
Consider first $|\rho|\ge\frac34\,n$. Then, $\sin\phi\le\frac{2}{3n}$,
thus $2n\,\phi\le2n\arcsin\frac{2}{3n}<\frac\pi2~~\forall n\ge1$.
Therefore $\real\p{1+\p{\frac1n\,\rho\,e^{i\phi}}^{2n}}>1$ and (\ref{br})
with $\ln C=0$.
Next, if $|\rho|<\frac34\,n$, then
$$
  \real\p{n^2-\sqrt[n]{n^{2n}+\p{\rho\,e^{i\phi}}^{2n}}}
        < n^2-n^2\sqrt[n]{1-\p{\tfrac34}^{2n}} < 1 ~.
$$
(The second function is decreasing for $n\ge2$.)~
Again, this implies (\ref{br}), with $\ln C=1$.
So we have shown that $\forall n\in\N^*,~\sigma_2(h_n)<3$, which
proves (\ref{de,mu}) for $h_n$.  With all that precedes, it is easy to
see that the given $\phi^n$ defines a net of $\{m,\pow\}$--mollifiers.

Now turn to $k_n\in\Sigma^\pow$. Estimating $x^\beta k_n(x)$
separately for $|x|\le 2$ and $|x|>2$ one can easily prove (\ref{de,mu}).
Once again, this allows to conclude that the
given $\phi^n$ defines a net of $\{m,\der\}$--mollifiers.
\end{proof}


\noindent The embedding of ultradistributions into the corresponding
weighted algebra of sequences  is realized through the first part of the
next theorem. Its second part deals with the representatives of
ultradifferentiable functions implying that  the multiplication
of regular enough elements within corresponding algebras is the
ordinary multiplication.


\begin{theorem}\label{thm35}
Assume $m>1$.
\begin{itemize}
\item[(i)]
Let $\psi\in\D^{(m)}$ (resp. $\psi\in\D^{\{m-\rho\}})$ with $\rho>0$
such that $m-\rho>1$) be compactly supported, and ${(\phi^n)}_n$ be a
net of $\{m,\pow\}$--mollifiers.  Then
\begin{eqnarray*}
  \psi * \phi_n - \psi &\in& \eN^{(p!^{m},p!^m)} ~,
  \qquad (\phi_n=n\,\phi^n(n\cdot))
\\\mbox{( resp. ~}
  \psi * \phi_n - \psi &\in& \eN^{\{p!^{m},p!^m\}}~\mbox)~.
\end{eqnarray*}

\item[(ii)]
Let $f\in\E^{\prime\,(m)}$ (resp.  $f\in\E^{\prime\,\{m\}}$) with
compact support; and $(\phi^n)_n$ a net of $\{m,\der\}$--mollifiers.
Then $f*\phi_n \in \E^{(p!^m,p!^{m-1})}_{exp}$,
(resp. $f*\phi_n \in \E^{\{p!^m,p!^{m-1}\}}_{exp}$).

If $(\phi^n)_{n}$ and $(\phi^{\prime n})_{n}$ are nets of
$\{m,\pow\}$--mollifiers,\comment{changed $der \to pow$}
then
\begin{eqnarray*}
  \forall\psi\in\D^{(m)}:
        \lr<>{f*\phi_n-f*\phi_n',\psi}&\in&\eN^{m}_0 ~,
\\\mbox{(~resp. ~}
  \forall\psi\in\D^{\{m-\rho\}}:
  \lr<>{f*\phi_n-f*\phi_n',\psi}&\in&\eN^{m}_0 ~)~.
\end{eqnarray*}
\end{itemize}
\end{theorem}

\begin{remark}
If $\psi\in\D^{(m)},\; m>1$, then
$(\psi)_n\in\E^{(p!^m,p!^{m'})}$ for all $m'>0$.  Fix a net of
$\{m,\pow\}$--mollifiers ${(\phi_n)}_n$.  The embedding
$\D^{(m)}\to\E^{(p!^m,p!^{m})}$ can be realized through $\psi\mapsto
{(\psi*\phi_n)}_n$ as well as through $\psi\to{(\psi)}_n$.  This is a
consequence of assertion (i). The similar conclusion follows for
$\D^{\{m-\rho\}}$. Thus, the product of $\varphi,\psi \in \D^{(m)}$
(resp. $\varphi,\psi \in \D^{\{m-\rho\}}$) is the usual one in
$\E^{(p!^m,p!^{m})}$ (resp. in $\E^{\{p!^m,p!^{m}\}}$).

\noindent
Assertion (ii) characterizes the embedding of elements in
$\E^{\prime\,(m)}$ (resp. $\E^{\prime\,\{m\}}$) into the
corresponding algebra by regularizations by
$\{m,\der\}$--mollifiers. Moreover, we have that the regularization of
elements in $\E^{\prime\,(m)}$ (resp. $\E^{\prime\,\{m\}}$) 
with $\{m,\pow\}$--mollifiers are weakly equal 
in the sense of ultracomplex numbers.

\noindent Note that $\D^{(m_1)} \hookrightarrow D^{\{m_1\}}\hookrightarrow D^{(m_2)}$,
$m_2>m_1>1$, where the left space is dense in the right one. This implies 
$\D^{\prime\,(m_2)} \hookrightarrow D^{\prime\,\{m_1\}}\hookrightarrow D^{\prime\,(m_1)}$. With
these relations theorem \ref{thm35} implies various embedding results
depending on the parameter $m>1$.
\noindent 
\end{remark}

\begin{proof}
\nb{(i)} Assume $\supp\psi\subset [-\mu,\mu]$.  Since $\psi*\phi_n-\psi=0$
for $|x|>\mu$, $n>n_0$, we assume in this proof $x\in [-\mu,\mu],~
n>n_0$.
First, we prove the assertion for the Beurling case; the Roumieu
case is treated in a similar way.

Let $s\in\N$. We have
\begin{eqnarray*}\lefteqn{
         (\psi\ast \phi_n-\psi)^{(s)}(x)
~=~     \INT_\R\p{\psi^{(s)}(x+t/n)-\psi^{(s)}(x)} \phi^n(t) dt
}\\&=& \INT_\R \left( \sum\limits_{p=0}^{N-1}
                \frac{ t^p }{ n^p\, p! } \psi^{(p+s)}(x)
                + \frac{ t^N }{ n^N N! } \psi^{(N+s)}(\xi)
                -  \psi^{(s)}(x)
        \right) \phi^n(t) dt ~,
\end{eqnarray*}
where $x\le\xi\le x+t/n$.  For $N=[n^{1/m}]+1$ as in the
definition of $\{m,\pow\}$--mollifiers,
\[
        (\psi*\phi_n-\psi)^{(s)}(x) = \INT_\R \frac{ t^N }{ n^N\,N! }
        \psi^{(N+s)}(\xi)\,\phi^n(t)dt ~.
\]
Let $b>1$ such that $\sigma^b(\phi^n)<\infty$. Then
\begin{eqnarray*}\lefteqn{
        \left| \frac{\nu^s}{{s!}^m} (\psi*\phi_n-\psi)^{(s)}(x) \right|
}\\&\le&
        \INT_\R \frac{1}{ {(N+s)!}^m }
        \left| \psi^{(N+s)}(\xi) \right|
        \frac{\nu^s {(N+s)!}^m}{n^N {s!}^m N!} t^N
        |\phi^n(t)| dt ~.
\end{eqnarray*}
We use ${N!}^m\le{(N^N)}^m$, $(N+s)!\le e^{N+s}\,N!\,s!$
~and~ $\frac{1}{n^N}\le\frac{2^N}{N^{Nm}}$~,
to get
\begin{eqnarray*} \lefteqn{
        \left|\frac{\nu^s}{s!^m}(\psi*\phi_n-\psi)^{(s)}(x)\right|
}&&\\&\le&
        \INT_\R \frac{(2e\,(\nu+b))^{N+s}}{(N+s)!^m}
        \left|\psi^{(N+s)}(\xi)\right| \frac{ N!^m }{ N^{mN} }
        \frac{ |t|^N }{ b^NN! } |\phi^n(t)| dt ~.
\end{eqnarray*}
Let $\ell>1$. Inserting $e^{-\ell N} e^{\ell N},$
with $ \nu_0 = 2\,\ell\, e\, (\nu+b)$, we have
\[
        \left| \frac{r^s}{s!^m} 
                \left(\psi*\phi_n-\psi\right)^{(s)}(x) \right|
        \le 2^{-\ell N} p^{m,\mu}_{\nu_0}(\psi)\,\sigma^b(\phi^n) ~.
\]
Now we use
$ e^{-\ell N} \sim e^{-\ell n^{1/m}}$ as $n\to\infty$.
This implies for every $\nu>0$ and $\ell>0$ there exist $C>0$ so that
\[
        \left| \frac{\nu^s}{{s!}^m}
                (\psi*\phi_n-\psi)^{(s)}(x) \right|
        \le C\,e^{-\ell\,n^{1/m}} ~.
\]
Taking the supremum over all $s$ and $x$, we obtain that
\[
  \ultra{ \psi * \phi_n-\psi }_{p^{m,\mu}_\nu,m} = 0    ~.
\]
\nb{Roumieu case:}
Let $d>1$  such that  $\sigma^d(\phi^n)<\infty $ and $h>0$
such that $p^{m-\rho,\mu}_{e^{m-\rho}h}(\psi)<\infty$.
We have, as above,
\begin{eqnarray*}       \lefteqn{
        \left| \frac{\nu^s}{s!^{m}} 
(\psi\ast \phi_n-\psi)^{(s)}(x) \right|
}&\\&\leq& \INT_{\R} 
        \frac{|\psi^{(N+s)}(\xi)|}{(N+s)!^{m-\rho} }
        \frac{\nu^s(N+s)!^{m-\rho}}{n^N s!^{m} N!} t^N
        |\phi^n(t)| d t ~.
\\&\leq& \INT_{\R} 
        \frac{(he^{m-\rho})^{N+s}|\psi^{(N+s)}(\xi)|}{(N+s)!^{m-\rho}} 
        \frac{ N!^{m}}{N^{Nm}}
        \frac{(h\nu)^ss!^{m-\rho}(dh)^N}{s!^m N!^\rho }
        \frac{|t|^N}{d^NN!}
        |\phi^n(t)| d t ~.
\end{eqnarray*}
Let $\ell>1$. Note
$$
        \sup\{\frac{(h\nu)^s s!^{m-\rho}}{s!^m}, s\in \N\}<\infty,~~
        \sup\{\frac{(dhe^\ell)^N}{N!^\rho}, N\in \N\}<\infty.
$$
As above we have, with suitable $C>0$,
(inserting $e^{-\ell N} e^{\ell N}$),
\[
        \lr||{ \frac{\nu^s}{{s!}^m}
                \p{\psi*\phi_n-\psi}^{(s)}(x) }
        \le C\,e^{-\ell N}\, p^{m-\rho,\mu}_{e^{m-\rho}h}(\psi)
                \,\sigma^d(\phi^n) ~.
\]
Again as above we finish the proof.

\nb{(ii)} We will give the proof in the Beurling case. The proof in the
Roumieu case is similar.
Recall~\cite{kom1}, if $f\in\E'{}^{(m)}$,
then there exists an ultradifferential operator of class $(m)$,
$P(D)=\sum_{k\in\N}
a_k\,D^k$, $\mu_0>0$
and continuous functions ${(F_k)}_{k\in\N}$,
with the property
$\supp F_k\subset[-\mu_0,\mu_0]$,
$
  \smash{\sup\limits_{ k\in\N, x\in\R }} |F_k(x)| \le M
$,
such that $f=\sum_{k\in\N}
a_k\,D^kF_k$.
%
This implies
$$
        \forall x\in\R: f*\phi_n(x) =
        \sum_{k=0}^\infty (-1)^k a_k\,n^k
        \INT_\R F_k(x+t/n)\,D^k\phi^n(t)dt ~,
$$
where ${(\phi_n)}_n$ is a net of $\{m,\der\}$--mollifiers 
such that $\sigma_b(\phi^n)<\infty$ and $a_k$, $k\in\N$
satisfy
$$
        \exists h,B>0:\forall k\in\N:|a_k|<Bh^k/k!^m ~.
$$
As in the part (i), we take $x\in[-\mu,\mu],\; \mu>\mu_0$ and $n>n_0$.
Let $\nu>1$ be given and $s\in\N$.  We have
\begin{eqnarray*}
{
  \frac{\nu^s}{s!^m} \left| f^{(s)} * \phi_n(x) \right|
}&=&
        \left|\sum_{k=0}^{\infty} (-1)^k a_k n^{k+s}
        \frac{\nu^s}{{s!}^m}
        \INT\limits_\R F_k(x+t/n)\,D^{k+s}\phi^n(t) dt\right|
\\&\le&
        \sum_{k=0}^{\infty}
        B\frac{\nu^s h^kn^{k+s}}{{k!}^m {s!}^m}
        \INT\limits_\R |F_k(x+ t/n)| \, |D^{k+s}\phi^n(t)| dt
\\&\le&
        \sum_{k=0}^{\infty}
        B\frac{(\nu h)^{s+k}n^{k+s}}{(k+s)!^{m}}
        \INT\limits_\R |F_k(x+t/n)|\,|D^{k+s}\phi^n(t)| dt
\\&\le&
        \sum_{k=0}^{\infty}\frac B{2^k}
        \frac{(2eb\nu h)^{s+k}n^{k+s}}{(k+s)!^{m-1}}
        \INT\limits_{\R}
        \frac{ | F_k(x+ \tfrac tn) | }{b^{k+s} \, (k+s)!}
        \left|D^{k+s}\phi^n(t)\right| d t
\\&\le&
        C\, e^{(2 e b \nu h n)^{1/(m-1)}} \sigma_b(\phi^n) ~.
\end{eqnarray*}
This proves that $f*\phi_n \in \E^{(p!^m,p!^{m-1})}_{exp}$.

Let us prove (for the Beurling case) that
$$
  \lr<>{ f,(\check\phi_n-\check\phi'_n)*\psi}\in\eN^{m}_0~.
$$
By continuity, we know that there exist $\mu\in\Ns$, $\nu>0$ and $C>0$
such that
\begin{eqnarray}\nonumber
  |\langle f,(\check\phi_n-\check\phi'_n) * \psi \rangle|
& \leq& C \, p^{\mu,m}_\nu((\check\phi_n-\check\phi'_n)* \psi)
\\
& \leq& C\,\lr[]{ p^{\mu,m}_\nu(\check{\phi}_n * \psi -\psi)+
  p^{\mu,m}_\nu(\check{\phi}'_n * \psi- \psi) } ~.\label{3}
\end{eqnarray}
By the first part of the theorem we have that
$$
  \psi*\phi_n-\psi ,~ \psi*\phi'_n-\psi \in\eN^{(p!^m,p!^m)} ~.
$$
This implies that for every $k>0$ there exists $C>0$ such that for
every $n\in\Ns$ both addents in (\ref{3}) are $\le C\,e^{-k\,n^{1/m}}$.
\end{proof}

\section{Generalized hyperfunctions on the circle}

For $\l>1$, let $\Omega_\l=\set{z\in\C\mid \frac1\l<|z|<\l}$ and
$\O_\l$ the Banach space of bounded holomorphic functions on
$\Omega_\l$. We denote by $\E(\T)$
(resp. $\A(\T):=\ind_{\l\to1}\O_\l$) the space of smooth
(resp. analytic) functions on the unit circle $\T=\{z\in\C\mid|z|=1\}$
and by $\E'(\T)$ (resp. $\B(\T)$) the corresponding space of
distributions (resp. hyperfunctions), cf.~\cite{MORI}.
For $f\in\A(\T)$, the coefficient $\widehat{T}(k)$ of $e_k(z)=z^k$ in
the Laurent expansion of $f$ is its $k$-th Fourier coefficient.
Complex numbers $c_k,\,k\in\Z$, are the Fourier coefficients of some
analytic function (resp. some hyperfunction) if and only if
$\ultra{\p{c_{\pm k}}_{k\in\N^*}}_{|\cdot|,1/k}<1$ (resp. $\le 1$).

\noindent Let $m\in [0,1)$ and $\nu>0$. We denote by $\A_{m,\nu}(\T)$
the set of functions $f\in\A(\T)$ such that $q^{m,\infty}_{\nu}(f)
:= \sup_{t\in\R,\alpha\in\N}\frac{|\tilde{f}^{(\alpha)}(t)|}
{\nu^{\alpha}{\alpha!}^m}<\infty$ where
$\tilde{f}(t)=f(e^{it}),~t\in\R$. We set
$$
        \A_m(\T)=\ind_{\nu \to \infty}\A _{m,\nu}(\T)
\text{~ and ~}
        \A_1(\T)=\ind_{m \to 1^-}\A_m(\T).
$$
Clearly $\A_1(\T)$ a subalgebra of $\A(\T)$ whose elements are
holomorphic in $\C^*$.


To prove the following theorem, we establish
\begin{lemma}\label{lemm aaa}
Let $m\in (0,1)$ and $\rho>e$. Let $\varphi$ denote the function
defined on $[\frac{1}{2},+\infty)$ by
$$
  \varphi(t)=\rho^{-t} \, t^{m(t+\frac{1}{2})} \, e^{-mt}.
$$
There exists a unique point $t_{\rho}\in [\frac{1}{2},+\infty)$ such
that $\inf_{t\ge\frac{1}{2}}\varphi(t)=\varphi(t_{\rho})$. Moreover
$\frac{1}{2} <\rho^{1/m}-t_{\rho}<\frac{1}{2}e^{\frac{1}{2t_{\rho}}}$.
It follows that $t_{\rho}\sim\rho^{1/m}-\frac{1}{2}\,\,(\rho\to
\infty)$ and
$$
  \sqrt{\rho}  e^{-\frac{m}{2}}  e^{-m\rho^{1/m}}
  < \varphi(t_{\rho})
  = \sqrt{\rho} e^{-m(t_{\rho}+\frac{1}{2}+\frac{1}{4t_{\rho}})}
  < \sqrt{\rho}e^{-m\rho^{1/m}}.
$$
\end{lemma}
\begin{proof}
Write $\psi(t) = -t\ln\rho+m \left( t+\frac{1}{2} \right) \ln t-mt$,
$ t \ge \frac{1}{2}$.
We find that the derivatives
$\psi'$ and $\psi''$ of $\psi$ are given by:
$$
  \psi'(t) = -\ln\rho+m\left(\ln t+\frac{1}{2t}\right);~
  \psi''(t) = \frac{m}{t}\left(1-\frac{1}{2t}\right).
$$
It follows that $\psi''(t)\ge 0$ for $t\ge\frac{1}{2}$ and then
$\psi(t)'$ vanishes at only one point $t_{\rho}$. Hence we have
$$
  \ln t_{\rho}+\frac{1}{2t_{\rho}}=\ln(\rho^{1/m}).
$$
It is seen that $\varphi$ increases on $[t_\rho,+\infty)$ and
$ \inf_{t\ge \frac{1}{2}}\varphi(t)= \varphi(t_{\rho})$.  \\
Let us compute $\varphi(t_{\rho})$. We have
$\psi(t_{\rho})=-t_{\rho}\ln\rho+m(t_{\rho}+\frac{1}{2})
\left(\ln\rho^{1/m}-\frac{1}{2t_{\rho}}\right)-mt_{\rho}$. We find
$\psi(t_{\rho})=\frac{1}{2}\ln\rho-\frac{m}{2}-\frac{m}{2t_{\rho}}-mt_{\rho}$
and $\varphi(t_{\rho})$ follows immediately.\\
{}From $\ln t_{\rho}+\frac{1}{2t_{\rho}}=\ln(\rho^{1/m})$ it is seen
that $t_{\rho} e^{\frac{1}{2t_{\rho}}}=\rho^{1/m}$, and then
$\rho^{1/m}-t_{\rho}=t_{\rho} (e^{\frac{1}{2t_{\rho}}}-1)$. Since
$x<e^x-1<xe^x$ for $x\ne0$, it follows that
$\frac{1}{2t_{\rho}}<e^{\frac{1}{2t_{\rho}}}-1
 < \frac{1}{t_{\rho}}e^{1/2t_{\rho}}$ and then
$\frac{1}{2} < \rho^{1/m}-t_{\rho} < \frac{1}{2}e^{1/2t_{\rho}}$,
because of $t_{\rho}\ge\frac{1}{2}$.\\
$t_{\rho}\ge\frac{1}{2}$ also implies that
$\rho^{1/m}-t_{\rho}<\frac{e}{2}$ showing that
$\lim_{\rho\to\infty}t_{\rho}=\infty$ and then
$\lim_{\rho\to\infty}e^{\frac{1}{t_{\rho}}}=1$.  From this, we find
$t_{\rho}\sim\rho^{1/m}-\frac{1}{2}$ $(\rho\to\infty)$.\\
It is seen that
$t_{\rho}+\frac{1}{2}+\frac{1}{4t_{\rho}}<\rho^{1/m}+\frac{1}{2}$;
from which it follows that
$\sqrt{\rho}e^{-\frac{m}{2}}e^{-m\rho^{1/m}}<\varphi(t_{\rho})$. Since
$t_{\rho}<\rho^{1/m}-\frac{1}{2}$ and $\varphi$ increases on
$[\frac{1}{2},+\infty[$, then
$\varphi(t_{\rho})<\varphi(\rho^{1/m}-\frac{1}{2})$. We seek for an
estimate to $\varphi(\rho^{1/m}-\frac{1}{2})$: \\
$\psi(\rho^{1/m}-\frac{1}{2})=-(\rho^{1/m}-\frac{1}{2})\ln\rho+m\rho^{1/m}
\ln(\rho^{1/m}-\frac{1}{2})-m(\rho^{1/m}-\frac{1}{2})$.  Writing
$\ln(\rho^{1/m}-\frac{1}{2})=\frac{1}{m}\ln\rho+\ln(1-\frac{1}{2\rho^{1/m}})$
and using $\ln(1-\frac{1}{2\rho^{1/m}})\le -\frac{1}{2\rho^{1/m}}$
yields $\psi(\rho^{1/m}-\frac{1}{2})\le
\frac{1}{2}\ln\rho-m\rho^{1/m}$ and consequently
$\varphi(\rho^{1/m}-\frac{1}{2})\le \sqrt{\rho}\,e^{-m\rho^{1/m}}$
proving the lemma.
\end{proof}

\begin{remark}\label{rema aaa}
It is easily seen that we also have
$\varphi(\rho^{1/m}+\frac{1}{2})\le \sqrt{\rho}e^{-m\rho^{1/m}}$.
\end{remark}

\def\ultrapm#1{\ultra{#1}^\pm}
\def\ultrapmk#1{\ultrapm{{(#1)}_k}}
\def\ultrapmmk#1{\ultrapmk{#1}_{(\cdot)^{-1/m}}}
\begin{theorem}\label{abe}
Let $f\in\A(\T)$ and $m\in(0,1)$.
\\
(i) If $f\in\A_{m,\nu}(\T)$ then:
$$
        \ultrapmmk{\hat{f}(k) }
        \le e^{-m/\nu^{1/m}}.
$$
Conversely if the above condition holds, then $f\in\A_{m,\nu'}(\T )$
for all $\nu'>\nu$.
\\
(ii) $f\in\A_{m}(\T )$ if and only if
$$
        \ultrapmmk{\hat{f}(k) } < 1.
$$
(iii) $f\in\A_{0,\nu}(\T)$ if and only if $\hat{f}(k)=0$ for $|k|>\nu$.\\
(iv) $f\in\A_{0}(\T )$ if and only if $(\hat{f}(k))_{k\in\Z}$ have
finite support.\\
(v) For all $f\in\A_1(\T)$ there exists $g\in\O(\C^*)$ such that $g|_\T=f$.
\end{theorem}

\begin{proof}
Let $f\in\A _{m,\nu}(\T)$ with $0<m<1$. For all $\alpha\in\N $,
$
  \tilde{f}^{(\alpha)}(t) =
  \sum_{p\in\Z} (ip)^{\alpha} \hat{f}(p) \, e^{ipt}
$. It follows that
$
  \INT_{-\pi}^{\pi}\tilde{f}^{\alpha}(t)e^{-ikt} dt
  = 2\pi(ik)^{\alpha}\hat{f}(k)
$,
and then consequently there is a positive constant $C_1$
such that $|k|^{\alpha}|\hat{f}(k)|\le C_1\nu^{\alpha} \alpha!^m$.
\\
Using Stirling's formula: $
  \alpha! = \alpha^{\alpha+1/2} \,
   e^{-\alpha}\sqrt{2\pi} \,
  (1+\varepsilon_{\alpha})
$,
$ \e_\alpha\searrow 0 $,
we find a positive constant $C_2$ such that:
$$
  \forall\alpha\in\N ^*,\,
  \forall k\in\Z,\,
  |k|^{\alpha}|\hat{f}(k)|
  \le C_2\nu^{\alpha} \alpha^{m(\alpha+1/2)}e^{-m\alpha}.
$$
It follows that:
$$
  \forall\alpha\in\N^*,\, \forall k\in\Z^*,\,
  |\hat{f}(k)|\le C_2 \left( \frac{\nu}{|k|} \right)^{\alpha}
    \, \alpha^{m(\alpha+1/2)} \, e^{-m\alpha}.
$$
Using the notations of Lemma \ref{lemm aaa} by taking
$\rho=\frac{|k|}{\nu}$ with $|k|>e\,\nu$, yields
$|\hat{f}(k)|\le C_2\,\varphi(t)$ for all $t\in\N^*$.
Following Remark \ref{rema aaa}, we have
$
  \varphi(\rho^{1/m}+\frac{1}{2}) \le
  \sqrt{\rho}e^{-m\rho^{1/m}}
$. Since $\varphi$ increases on $
  [\rho^{1/m}-\frac{1}{2},\rho^{1/m}+\frac{1}{2}]
$ which contains a positive integer $\alpha_{\rho}$,
then $
  |\hat{f}(k)| \le C_2\,\varphi(\alpha_{\rho})
  \le C_2\sqrt{\rho} \, e^{-m\rho^{1/m}}
$ for $|k|>e\nu$, that is 
$
  |\hat{f}(k)|\le C_2 \left( \frac{|k|}{\nu} \right)^{\frac{1}{2}}
    \, e^{-m/\nu^{1/m}}
$
for $|k|>e\nu$ from which inequality of (i) follows.\\
Conversely assume that $f$ satisfies the condition of (i).
Let $\nu'>\nu$. Choose $\nu''$ such that $\nu'>\nu''>\nu$ and set
$\beta'=m/(\nu')^{1/m}$, \, $\beta''=m/(\nu'')^{1/m}$. 
Since $e^{-\beta''}>e^{-m/\nu^{1/m}}$, there exists a positive constant 
$C>|\hat{f}(0)|$ such that $|\hat{f}(k)|\le Ce^{-\beta''|k|^{1/m}}$ for every
$k\in\Z$. 
For every $\alpha\in\N $, we have
$
  \tilde{f}^{(\alpha)}(t) =
  \sum_{k\in\Z}(ik)^{\alpha}\hat{f}(k) \, e^{ikt}
$.
Then, using the last iequality we find
$$
    \|\tilde{f}^{(\alpha)}\|_{\infty}
    \le C \, \left(\sum_{k\in\Z}
        e^{-(\beta''-\beta')|k|^{1/m}}\right)
    \sup_{k\in\Z } |k|^{\alpha} \,
    e^{-\beta'|k|^{1/m}}.
$$
Let $\phi(t)=t^{\alpha}e^{-\beta't^{1/m}};\, t\ge 0$.
A simple study of $\phi$ shows that
$
  \sup_{t\ge0}\phi(t) = \phi(\nu'\alpha^m)
  = (\nu')^{\alpha}\alpha^{m\alpha}
  \, e^{-m\alpha}
$.
Using Stirling's formula, we get a positive constant $C_1$
such that for all $ \alpha\in\N $,
$
  \left\| \tilde{f}^{(\alpha)} \right\|_\infty
  \le C_1 \, \nu'{^\alpha} \alpha!^m
$,
showing that $f\in\A_{m,\nu'}(\T)$ and proving (i).

Let $f\in\A_{m}(\T )$. Then $f\in\A_{m,\nu}(\T )$
for some $\nu>0$ and the inequality follows from (i) and $e^{-m/\nu^{1/m}}<1$.
Conversely if $\ultrapmmk{\hat{f}(k)}<1$,
there exists $\nu>0$ such that $\ultrapmmk{\hat{f}(k)}
\le e^{-m/\nu^{1/m}}$. From (i), it follows that $f\in\A_{m,\nu'}
(\T)$ for $\nu'>\nu$. Hence $f\in\A_{m}(\T)$ proving (ii).

Let $f\in\A_{0,\nu}(\T)$. The previous shows that
there exists $C_1>0$ such that $
  |k|^{\alpha}|\hat{f}(k)| \le C_1\nu^{\alpha}
$. It follows that $
  |\hat{f}(k)| \le
    C_1 \left(\frac{\nu}{|k|}\right)^\alpha
$ for all $k\in\Z^*$ and all $\alpha\in\N $.
If $|k|>\nu$, then $\frac{\nu}{|k|}<1$, and making $\alpha\to\infty$
yields $\hat{f}(k)=0$.\\
Conversely, assume that $\hat{f}(k)=0$ for $|k|>\nu$. Then we have
$
  \forall z\in\C^*,\, f(z)=\sum_{|k|\le \nu}\hat{f}(k) \, z^k.
$
It follows that for all $\alpha\in\N $,
$ \left\| \tilde{f}^{(\alpha}) \right\|_{\infty} \le
\left( \sum_{|k|\le \nu} \left|\hat{f}(k) \right| \right) \,
\nu^{\alpha} $, that is $f\in\A_{0,\nu}(\T)$ proving (iii).\\
Claim (iv) follows from (iii) straightforwardly.\\
Claims (ii) and (iv) show that for $f\in\A_1(\T)$ the series 
$\sum_{k\in\Z}\hat{f}(k)\,z^k$ converges absolutely for any $z\in\C^*$,
proving (v).
\end{proof}

Now let $r={(r_n)}_n$ with $r_n>0$ and $r_n \searrow 0$.  For
$n\in\N$, we set $\psi_n(z)=\sum_{|k|\le 1/r_n}z^k$.  We have
$\psi_n*\psi_n=\psi_n$ and $\lim_{n\to\infty}\psi_n=\delta$
in $\E'(\T)$. If $H\in\B(\T)$, $H*\psi_n=\sum_{|k|\le
1/r_n} \hat{H}(k)z^k$ (where $S*T=z\mapsto
\sum_{k\in\Z}\widehat{S}(k)\,\widehat{T}(k)\,z^k$) and consequently
$\lim_{n\to\infty} H*\psi_n=H$ in $\B(\T)$.\\
For $f\in\O_\l$ we set
$
        q^\l(f) =
        \left\|\,f\,\right\|_{L^\infty(\Omega_\l)}
$
and 
$
\hat q^\l(f)=\sup\limits_{k\in\Z}\l^{|k|}|\hat f(k)|.
$
For a sequence $f={(f_n)}_n\in(\O_\l)^\N$,  we set
$
        \ultra{f}_{q^\l,r} :=
        \limsup\limits_{n\to\infty} {q^\l(f_n)}^{r_n}
$
and 
$
        \ultra{f}_{\hat q^\l,r} :=
        \limsup\limits_{n\to\infty} {q^\l(f_n)}^{r_n}.
$
We  define 
\begin{align*}
        \overrightarrow\F_{q,r} &= \set{ f\in\A_1(\T)^\N \mid
        \exists\l>1: \ultra{f}_{q^\l,r} < \infty } ~,
\\
        \overrightarrow\K_{q,r} &= \set{ f\in\A_1(\T)^\N \mid
        \exists\l>1: \ultra{f}_{q^\l,r} = 0 } ~.
\\
        \overrightarrow\F_{\hat{q},r} &= \set{ f\in\A_1(\T)^\N \mid
        \exists\l>1: \ultra{f}_{\hat{q}^\l,r} < \infty } ~,
\\
        \overrightarrow\K_{\hat{q},r} &= \set{ f\in\A_1(\T)^\N \mid
        \exists\l>1: \ultra{f}_{\hat{q}^\l,r} = 0 } ~.
\end{align*}

Then we have
\begin{proposition}\label{aba}
Let $\l>1$ and $f=\left(f_n\right)_n\in\O_\l^\N$. Then we have for 
any $\mu\in(1,\l)$
$$
        \ultra f_{q^\mu,r}\le \ultra f_{\hat q^\l,r}
\le \ultra f_{q^\l,r}.
$$
Then consequently $\overrightarrow\F_{q,r}=\overrightarrow\F_{\hat q,r}$ 
and $\overrightarrow\K_{q,r}=\overrightarrow\K_{\hat q,r}$ 
\end{proposition}

\PROOF{Let $\l>\mu>1$ and $f=\left(f_n \right)_n \in\O_\l^\N$.
For every $k\in\Z$, $\l^{|k|}|\hat{f}_n(k)|\le\hat q^\l(f_n)$ and then 
$|\hat{f}_n(k)|\le\hat q^\l(f_n)\l^{-|k|}$. Using this inequality, 
we find from  $|f_n(z)|\le 
\sum_{k\in\Z}|\hat{f}_n(k)|\,|z|^k$ that $|f_n(z)|\le \hat q^\l(f_n)
\sum_{k\in\Z}\left(\frac{\mu}{\l}\right)^k$. It follows that there exists a positive 
constant $C(\l,\mu)$ such that $C(\l,\mu)q^\mu(f_n)\le \hat q^\l(f_n)$.\\
From Cauchy's formula $\l^{|k|}|\hat{f}_n(k)|\le q^\l(f_n)$ and then 
$\hat q^\l(f_n)\le q^\l(f_n)$. Finaly we obtain
$$
        \ultra f_{q^\mu,r}\le \ultra f_{\hat q^\l,r}
\le \ultra f_{q^\l,r},
$$
and then it follows srtraightfowardly that $\overrightarrow\F_{q,r}=
\overrightarrow\F_{\hat q,r}$ and $\overrightarrow\K_{q,r}=
\overrightarrow\K_{\hat q,r}$.}
\begin{definition}
Let $\X_r(\T) =\overrightarrow\F_{q,r}$ and $\eN_r(\T) =
\overrightarrow\K_{q,r}$. The algebra of generalized hyperfunctions on 
$\T$ is $\G_{H,r}=\X_r(\T)/\eN_r(\T)$.
\end{definition}
We have an embbeding of  $\B(\T)$ in $\G_{H,r}(\T)$ which preserves 
the usual multiplication of elements in $\A_1(\T):$
\begin{theorem} Let
$$
\begin{array}[t]{lrcl}  
        {\bf\bar{i}}: & \B(\T) & \to & \G_{H,r}(\T) \\
                 & H & \mapsto & \left[ {(H*\psi_n)}_n \right]
\end{array}
\text{ and ~ }
\begin{array}[t]{lrcl}  
        {\bf\bar{i}_0}: &\A_1(\T) & \to & \G_{H,r}(\T) \\
                & f & \mapsto & \left[ {(f)}_n \right] \end{array}.
$$
Then, $\bf\bar{i}$ is a linear embedding and $\bf\bar{i}_0$ is a
one to one morphism of algebras such that
${\bf\bar i}|_{\A_1(\T)}={\bf\bar i_0}$.
\end{theorem}

\begin{proof}
The claim on $\bf\bar{i}_0$ is
easy to prove. Let us focus on the properties of the
first part related to $\bf\bar{i}$. The linearity of
$\bf\bar{i}$ is quite obvious. Let $H\in\B(\T)$ and set
$h=(h_n)_n$ with $h_n=H*\psi_n$. 
%
%
{}From Theorem~\ref{abe},(v), we have $h\in\X(\T)$.

Now take $\l>1$. From the property of the Fourier
coefficients of $H$, there exists $C>0$ such that $|\hat{H}(k)|\le
C\,\l^{|k|}$ for all $k\in\Z$. It follows that 
$\l^{|k|}|\hat{h}_n(k)|\le {C}\l^{2/r_n}$ showing that
$\ultra{h}_{\hat{q}^{\l},r}\le\l^2$.
By Proposition~\ref{aba}, $h\in\X_r(\T)$.
It is sufficient to consider restrictions to the spaces $\A_m(\T)$
with $0<m<1$.  Let $f\in\A_m(\T)$ with $0<m<1$. There is $\l>1$
such that $f(z)=\sum_{k\in\Z} \hat{f}(k)\,z^k$ for
$1/\l\le|z|\le\l$.
Then we have ${\bf\bar{i}_0}(f)-{\bf\bar{i}}(f)=[f_n]$ where
$f_n=f-f*\psi_n$, that is
$f_n(z)=\sum_{|k|>1/r_n}\hat{f}(k)\,z^k$. Then we have
${(f_n)}_n\in\O_\l$.\\
We claim that ${(f_n)}_n\in\eN_r(\T)$.
{}From Theorem~\ref{abe}, there exist $p\in(0,1)$ and $C>0$ such that 
every $k\in\Z,\,\,|\hat{f}(k)|\le C p^{|k|^{1/m}}$. For $|k|>r_n^{-1}$,  
writing $p^{|k|^{1/m}}\le p^{\frac{1}{2}|k|^{1/m}}
p^{\frac{1}{2}r_n^{-1/m}}$, we find 
$
\left(\l^{|k|}|\hat{f}_n(k)|\right)^{r_n}\le \left(C\l^{|k|}
p^{\frac{1}{2}|k|^{1/m}}\right)^{r_n}p^{\frac{1}{2}r_n^{(m-1)/m}}
$
Since $C\l^{|k|}p^{\frac{1}{2}|k|^{1/m}}$ is bounded with respect to $k$, 
because of $1/m>1$ and $p\in (0,1)$, it follows that 
$
  \ultra{f}_{\hat{q}^{\l},r}=0 
$,
proving our claim.
\end{proof}


\bigskip\parindent0pt\parskip1ex\flushleft

Antoine Delcroix: IUFM de la Guadeloupe,\\
Morne Ferret, BP 399, 97159 Pointe \`a Pitre cedex (Guadeloupe, F.W.I.)\\
Tel.: 00590 590 21 36 21, Fax : 00590 590 82 51 11,\\
E-mail: {\tt Antoine.Delcroix@univ-ag.fr}

Maximilian Hasler: Universit\'e des Antilles et de la Guyane, D.S.I.,\\
BP 7209,
97275 Schoelcher cedex (Martinique, F.W.I.)\\
Tel.: 00596 596 72 73 55, Fax : 00596 596 72 73 62,\\
e-mail : {\tt Maximilian.Hasler@martinique.univ-ag.fr}

Stevan Pilipovi\'c: University of Novi Sad, Inst. of Mathematics, \\
Trg D. Obradovi\'ca 4, 21000 Novi Sad (Yougoslavia)\\
Tel.: 00381 21 58 136, Fax: 00381 21 350 458,\\
e-mail: {\tt pilipovic@im.ns.ac.yu}

Vincent Valmorin: Universit\'e des Antilles et de la Guyane, D.M.I.,\\
Campus de Fouillole, 
97159 Pointe \`a Pitre cedex (Guadeloupe, F.W.I.)\\
Tel.: 00590 590 93 86 96, Fax: 00590 590 93 86 98,\\
e-mail: {\tt Vincent.Valmorin@univ-ag.fr}


\end{document}